\documentclass[reqno]{amsart}

\usepackage{upref}
\usepackage{amsfonts, amsmath}
\usepackage{epsfig}
\usepackage{graphics, latexsym,amssymb}
\def\S3{{\mathbb S^3}}
\def\R3{{\mathbb R^3}}

\newtheorem{theorem}{Theorem}[section]
\newtheorem{proposition}{Proposition}[section]
\newtheorem{lemma}{Lemma}[section]
\newtheorem{definition}{Definition}[section]
\newtheorem{remark}{Remark}[section]
\newtheorem{corollary}{Corollary}[section]

\begin{document}

\title{A proof of the Lawson conjecture for minimal tori embedded in $\S3$}

\author{Fernando A. A. Pimentel}

\begin{abstract}
A peculiarity of the geometry of the euclidean 3-sphere $\S3$ is that it allows for the existence of compact without boundary minimally immersed surfaces. Despite a wealthy of examples of such surfaces, the only known tori minimally embedded in $\S3$ are the ones congruent to the Clifford torus. In 1970 Lawson conjectured that the Clifford torus is, up to congruences, the only torus minimally embedded in $\S3$. We prove here Lawson conjecture to be true.
Two results are instrumental to this work, namely, a characterization  of the Clifford torus in terms of its first eingenfunctions (\cite{MR}) and the assumption of a ``two-piece  property" to these tori: every equator  divides a torus minimally embedded  in $\S3$ in exactly two connected components (\cite{Rs}).

\end{abstract}

\maketitle

\allowdisplaybreaks

\section{Introduction}

A special feature of the 
the euclidean
3-sphere $\S3=\{x\in \mathbb R^4:|x|=1\}$ on which it differs from the euclidean space $\R3$ is that its topology does not obstruct the existence of compact without boundary minimally immersed surfaces. In fact, a plethora of examples of such surfaces are known to exist. Thus it was shown in \cite{hl} by Hsiang-Lawson the existence of infinitely many immersed tori. On the other hand, examples of embedded minimal surfaces of arbitrary genus are constructed in \cite{L1} by Lawson, with additional examples of such surfaces being suplemented by Karcher-Pinkall-Sterling in \cite{kps}.

As regards embedded minimal surfaces, any impression of superabundance in examples sugested by the papers cited above evanesces as soon as one realizes that all minimally embedded n-tori hitherto known must have plenty of symmetries, that are inherited from their construction process. For instance, the only known minimally embedded torus is, up to congruences, the Clifford torus (see Def. \ref{dctorus2}), which has all the symmetries compatible with its topology.
     
Some results indicate that a minimally embedded torus must indeed  be very symmetrical. Thus Lawson \cite{L3} proved that any embedded minimal torus must be unknotted. In 1995, Ros \cite{Rs} extended the restrictions on minimally embedded tori by proving that they have a special kind of {\it two-piece property} (see also Def. \ref{cl0}):
\begin{quotation} 
\it If $M$ is an embedded compact minimal surface in $\S3$, then
every totally geodesic equator of $\S3$ divides $M$ in exactly two connected components.
\end{quotation}
This last result will be used extensively in this work. 

It is well known that the euclidean coordinates $x_1,\dots,x_4$ of a minimal immersion $M$ in $\S3$
are eigenfunctions of the Laplacian on $M$ (see section \ref{sl}). With respect to this, we state a characterization of the Clifford torus due to Montiel-Ros \cite{MR} that shall be useful to us as well:
\begin{quotation} {\it The only minimal torus immersed into $\S3$ by its first eigenfunctions is the Clifford torus.}
\end{quotation}
This remarkable result  relates the Yau conjecture, which asserts that any minimal embedding of a compact surface
is by the first eigenfunctions and the Lawson conjecture:
\begin{quotation} CONJECTURE (Lawson, \cite{L3}): {\it a torus minimally embedded in $\S3$ is congruent to the Clifford torus}.
\end{quotation}
Statements of these conjecures may be found in \cite{y} (both conjectures) and \cite{L3} (original statement of Lawson conjecture). Our aim here is to produce a proof of Lawson conjecture (Theorem \ref{tl1}).  

Our approach to this problem is twofold, based on the results cited above in \cite{Rs} and \cite{MR}. Along Section \ref{Sptpp} we use the {\it two-piece property} to classify the intersection of an equator $S(v)$ with a torus $M$
with the {\it two-piece property} and strictly negative Gauss-Kronecker curvature (see Section \ref{Sct} and Def. \ref{cl0} for notation and basic definitions) in four types: 
\medskip

\noindent
{\bf Proposition \ref{pclass}}
{\it Let $M$ be a $C^{2,\alpha}$torus embedded in $\S3$ that has  strictly negative
Gauss-Kronecker curvature. Then $M$ has the  {\it two-piece property} if and only if for every equator $S(v)$ the intersection $S(v)\cap M$  
is classified among the four types listed below:
\begin{enumerate}
\item an embedded closed curve,
\item two disjoint  embedded closed curves,
\item two embedded closed curves bounding disjoint discs in $S(v)$ with only one  point in common,
\item two embedded closed curves bounding disjoint discs in $S(v)$ with exactly two points in common.
\end{enumerate}
where a curve in a intersection of type 1 is nullhomotopic in $M$ and the pair of curves in intersections of
types 2, 3 or 4 are homotopic in $M$ to a same generator of $\pi_1(M)$ for any choice of the pair of curves in type 4 intersections (see Remark \ref{rchico}).}
\medskip

\noindent
Diagrams of these intersections are sketched in Appendix \ref{afig}.

Then we show that 
for each torus $M$ as above there exists an equator $S(v)$  such that the components of $M\backslash S(v)$ are two tubes with common border disjoint smooth closed curves in $S(v)$:
\medskip

\noindent
{\bf Lemma \ref{ltype}} 
{\it There exists an equator $S(v)$ such that $S(v)\cap M$ is of type 2.}\medskip

These results establish a kind of ``topological simmetry'' strong enough to provide us the geometric support needed in our next step: granted the existence
of a minimal torus not congruent to the Clifford torus, construct a sequence of minimal tori converging in $C^{2,\alpha}$ to the Clifford torus. The {\it two-piece property} will be used to control the convergence process from which the sequence above shall be built (see Section \ref{scon}). We then prove Lemma \ref{lfinal}
and its corollary:\medskip

\noindent
{\bf Corollary \ref{cfinal}} {\it Let $0<\alpha<1$. If there exists a torus $M$ minimally embedded in $\S3$  that is not congruent to the Clifford torus $T$ then there also exists a
sequence $(X_k)$  of $C^{2,\alpha}$ diffeomorphisms of $\S3$ canonically extended to $A_2$ such that each $M_k=X_k(T)$
is a torus minimally embedded in $\S3$ noncongruent to the Clifford torus 
and $\lim_{k\to\infty}||X_k-I||_{C^{2,\alpha}}= 0$.}\medskip

\noindent
Above it is mentioned  the annulus $A_2=\{x\in \mathbb R^4 \hspace{.06in}|\hspace{.06in}1/2<|x|<2\}$ in connection with a canonical extension of maps of $\S3$ to maps of $A_2$, defined in Section \ref{scon}. This seemingly extraneous apparatus is in fact pertinent. The extension of maps 
from the ambient space $\S3$ to $A_2$ will somewhat simplify our approach due to the vector space structure of $\mathbb R^4$ thus incorporated. This will enable us to make use of ordinary Holder spaces of functions on regions in a euclidean space.

A crucial point in the proof of Lemma \ref{lfinal} and its Corollary \ref{cfinal} is to ensure that the tori $M_k$ are not congruent to the Clifford torus. With respect to this, we remark the role of Lemma \ref{ltype3}. It asserts that some closed curve in a type 2 intersection given by Lemma \ref{ltype} stated above is not congruent to any closed
curve in the Clifford torus $T$, assumed $M$ minimal and noncongruent to $T$. This fact will be used to distinguish the tori $M_k$ from the Clifford torus. 

On the other hand, from the supracited characterization of the Clifford torus in \cite{MR}, it will be proven in Section \ref{sl} that the Clifford torus is isolated in $C^{2,\alpha}$. In order to do that, it will be shown that whenever $(M_k)$ is a sequence of embedded tori noncongruent to the Clifford torus $T$  converging to $T$ in $C^{2,\alpha}$ then their first eigenvalues will converge to two, the first eigenvalue of $T$. But the first eingenfunctions of $M_k$ will not converge to the first eigenfunctions of  $T$ as well. From these remarks a  contradiction will arise and the result below will follow:\medskip

\noindent
{\bf Lemma \ref{ltiso}} {\it Let $0<\alpha<1$. Then there does not exist any sequence $(X_k)$ of diffeomorphisms of 
$\S3$ canonically extended to $A_2$  such that each $M_k=X_k(T)$ is a torus minimally embedded in $\S3$ noncongruent to the Clifford torus $T$ and $\lim_{k\to \infty}||X_k-I||_{C^{2,\alpha}}=0$.} \medskip

The main result in this paper (Theorem \ref{tl1}) is an immediate consequence of Corollary \ref{cfinal} and Lemma \ref{ltiso}. \bigskip

Part of this work was in the Doctoral thesis of the  author at
UFC, Brazil. I would like to thank A.
Gerv{\'a}sio Colares for his orientation and continuous encouragement. I would also like to thank
H. Rosenberg  and A. Ros for some very valuable discussions. I am specially grateful to Abd\^ enago A. de
Barros that introduced me to the Lawson conjecture and to the results in \cite{Rs} and \cite{MR}. Finally, I would like to express my gratitude to
the Universit\' e Paris XII and the Minist\` ere de la Recherche, France, for the post-doc position that enabled me to  discuss
my issues with some of the afore mentioned professors.

\section{Preliminaries}\label{Sct}

The elements of spherical geometry that we will need are sketched
here and found in full detail in \cite{bg}. 

Let $\S3\subset \mathbb R^4$ be the set of unit vectors of $\mathbb R^4$.
A totally geodesic equator, hereafter called equator, is a totally geodesic big 2-sphere of $\S3$. Let $\langle\hspace{.02in},\hspace{.02in}\rangle$ be the inner product of $\mathbb R^4$.
Following \cite{Rs}, we let  correspond to any unit vector $v\in\mathbb R^4$   the equator $S(v)=\{p\in \S3 \hspace{.02in}|\hspace{.02in}\langle v,p\rangle=0\}$ and the half spheres $H_+(v)=\{p\in \S3\hspace{.02in} |\hspace{.02in}\langle v,p\rangle >0\}$
and $H_-(v)=\{p\in \S3\hspace{.02in} |\hspace{.02in}\langle v,p\rangle <0\}$. 

For ``distance between points" it is meant here the intrinsic distance of $\S3$.
A circle $C$ in an equator $S(v)$ of $\S3$ is a closed curve in $S(v)$ whose points
equidist from a center, that is, a distinguished point in $S(v)$.\label {center}
According to this definition, the circle $C$ has two centers $p_1$, $p_2$ in the
equator $S(v)$, which are antipodal points (i.e., $p_1=-p_2$) as it is
easily seen.  The curvature of  $C$ is
given by $|\hbox{cotg}(r)|$, where $r$ is the distance from $C$ to
any of its centers (see in \cite{Sp}, vol 3, the definition of curvature of a curve). If $C$ has positive curvature then we can
define the normal vector field along $C\subset S(v)$ as the unit vector
field along $C$, normal to $C$ and tangent to $S(v)$, that points
towards the connected component of $S(v)\backslash C$ containing
the center closer to $C$.

 When a circle $C$ has zero curvature it is called
a geodesic or great circle. We can perform a rotation of an equator in $\S3$ around
any geodesic in it, so
a geodesic is contained in infinitely many equators with
a pair of centers in each equator.

Tori dealt with in this paper are  $C^{2,\alpha}$ maps $X:
S^1\times S^1\rightarrow \S3$, immersions (when $dX$ is injective) and embeddings (when it is further assumed that $X$ is injective).
Defined along these
tori are the principal curvatures, denoted by $k_1$ and
$k_2$, the mean curvature $H=(k_1+k_2)/2$ and the Gauss-Kronecker curvature $S=k_1k_2$. We recall that the  intrinsic or sectional
curvature of surfaces  in $\S3$ is given by $K=1+S$. Hence  the Clifford torus, with principal curvatures $k_1=1$ and $ k_2=-1$ for a suitably choice of a normal vector field is a minimal ($H\equiv 0$) flat $(K\equiv 0)$ torus as we will verify below. 
 We refer to \cite{Sp}, vols. 3 and 4, for the
definitions and a throughout analysis of the curvature functions along immersed surfaces.

\begin{definition}\label{dctorus2} We identify $\mathbb R^4$ with $\mathbb C^2$
in the usual manner. The  Clifford torus is the compact surface
in $\S3$ given by
$$\{(u,v)\in \mathbb C^2:\hspace{.03in} |u|=|v|=\sqrt 2/2 \}.$$
\end{definition}

The orbits of the action of the groups
$O^+(\mathbb C)\times\hbox{id}_{\mathbb C}$ and
$\hbox{id}_{\mathbb C} \times O^+(\mathbb C)$ in the Clifford
torus( see \cite{bg}, 18.8.6) are two families of circles parameterized by $u\rightarrow (u,v_0)$ and 
$v\rightarrow(u_0,v)$ for $|u|=|v|=|u_0|=|v_0|=\sqrt 2/2$. These circles have curvature equal to 1 and circles in different families are orthogonal at their only common point. Moreover, a simmetry argument implies that a vector field normal to any of these circles is normal to the Clifford torus too. Thus these families of curves are the two systems of lines of curvature of the Clifford torus.

\section {the two-piece property} \label{Sptpp}

Below it is presented a property concerning equators and subsets of $\S3$ already mentioned in the Introduction.

\begin{definition} \label {cl0}
We say that a set $S\subset \S3$ has the
{\it two-piece property}   when $S\backslash S(v)$
has exactly 2 connected components for each equator $S(v)$.
\end{definition}

The Clifford torus is itself an instance of surface with the {\it two-piece property}. In order to verify this either apply Ros' Theorem (see Introduction) or follow this straightforward reasoning: cyclides of Dupin (i.e., tori whose lines of curvature in the two systems are circles, e.g., any tori of revolution) are known to be divided in two connected components by any equator transversal to them at some point (see \cite{Ba}).  The Clifford torus is a cyclide of Dupin (see Section \ref{Sct}) transversal to any equator (since it has strictly negative Gauss-Kronecker  curvature (see Sect. \ref{Sct}), it is not contained in any halfspace $H_+(v)$ or $H_-(v)$; or one could move the equator $S(v)$ until it touches the Clifford torus by one of its sides).

 Nevertheless, being a cyclide of Dupin is so strong a constraint for a tori to have the
{\it two-piece property}. Hereafter we will obtain weaker conditions that imply the {\it two-piece property}  for tori in $\S3$.

Firstly, we remark that tori with the {\it two-piece property}\ can not have points of positive Gauss-Kronecker curvature: consider an equator tangent to a torus at a point of positive Gauss-Kronecker curvature. A suitable motion of the equator shall increase 
the number of connected components of the torus bordered by the equator. 

Besides,  minimally immersed tori do not have umbilical points (i.e., where $k_1=k_2=0$), see Lemma 1.4 and the ensuing comments in \cite{L1}. So these tori have strictly negative Gauss-Kronecker curvature too. Consequently, we will restrict our analisys to tori with strictly negative Gauss-Kronecker curvature. 

Now we observe that it can be proved that
the properties of tori that concern us here, i.e., embeddedness, negativity of the Gauss-Kronecker curvature and the {\it two-piece property}, rest undisturbed after small perturbations. Thus

\begin{lemma}\label{ltppopen}The set of embedded $C^{2,\alpha}$ tori with {\it two-piece property} and strictly negative Gauss-Kronecker curvature is open in the space of 
$C^{2,\alpha}$ tori endowed with the $C^{2,\alpha}$ topology.  
\end{lemma}

\begin{remark}\label{rtopos} A general definition of the $C^{2,\alpha}$ topology for maps from $M$ to  $N$, where $M$, $N$ are manifolds, can be found in \cite{Hi}, section 2.1.
\end{remark}

The following Lemmas  classify the curves in the intersection  of
equators and tori with the {\it two-piece property} and strictly negative
Gauss-Kronecker curvature embedded in $\S3$ in terms of their number, shape and homotopy class (we refer to \cite{Ma} for the nomenclature and results assumed hereon).

The terminology ``meridians'' and ``longitudes'' is well established and refers to the two systems of canonical generators of the fundamental group $\pi_1(\overline M)$ of a torus $\overline M$ embedded in $\mathbb R^3$. Informally, they are the curves that run once in the proper sense around the torus, thus linking any 
homotopically nontrivial curve in the bounded component of $\mathbb R^3 \backslash \overline M$ (meridians)
or in its unbounded component (longitudes) (see \cite{Ad}, p.108).

Analogously, there exists a natural system of 
generators of the fundamental group $\pi_1(M)$ of a  torus $M$ in $\S3$. They may be called  meridians and longitudes as its counterparts on  tori embedded in $\mathbb R^3$, after a choice of a normal vector field along $M$. We  define them as the images of the lines of curvature of the Clifford  torus $T$ by some diffeomorphism $X$ of $\S3$ such that $X(T)=M$ that preserves the chosen orientation. 
{\it For the sake of simplicity we shall adopt  ``generator'' for ``canonical generator''}.

For a disc it is meant any subset homeomorphic to $\{x\in \mathbb R^2|\hspace{.03in} |x|<1\}$ and a tube is 
a disc minus one point. {\it In what follows $M$ is an embedded $C^{2,\alpha}$ torus, $0<\alpha\leq 1$, with strictly negative Gauss-Kronecker curvature and the {\it two-piece property}}. 

\begin{lemma}\label{lclass1}
 If $S(v)$ is an equator that is not tangent to $M$ at any point then $S(v)\cap M$ either is an only embedded closed curve bounding a disc in $M$ or  is the union of two disjoint embedded closed curves
homotopic in $M$ to a same generator of $\pi_1(M)$.
\end{lemma}

\noindent {\it Proof.} Obviously $S(v)\cap M$ is not empty. In fact, it is the union of disjoint closed embedded regular curves since $M$ is an embedded compact smooth surface, and $M$ is strictly transversal to $S(v)$. 

Let $M^+$, $M^-$ be  the closure of the connected components of $M\backslash S(v)$, with $M^+\backslash S(V)\subset H^+(v)$ and $M^-\backslash S(v) \subset H^-(v)$. Thus, by the {\it two-piece property} and the embeddedness assumed for $M$, both $M^+$ and $M^-$ are compact connected orientable surfaces with boundary.
It is well known that $M^+$ and $M^-$ are then embedded orientable compact surfaces in $\S3$ without a finite number of disjoint  discs 
whose borders, i.e., $S(v)\cap M^+$ and $S(v)\cap M^-$, are the embedded closed curves in the intersection  $S(v)\cap M$ (see prob. 4, p. 144, in \cite {Sg}). We classify below $S(v)\cap M$ according to the topology of $M^+$.
\bigskip

\noindent
(1) $M^+$ is an sphere without an only disc. Then $S(v)\cap M$ contains an only embedded closed curve bounding a disc in $M$, that is,   
$S(v)\cap M^+$; 
\medskip\newline
(2) $M^+$ is an sphere without exactly two disjoint discs. Then $M^-$ must be itself an sphere without two discs,  otherwise $M$ is not a torus. Thus the intersection of $M^+$ with $S(v)$, i.e., $S(v)\cap M$, is the disjoint union of two curves homotopic to a same generator of $\pi_1(M)$;\medskip\newline
(3) $M^+$ is an sphere without $n\geq 3$ discs. Then at least $M^-$ is  an sphere without $n$ discs and $M$ is an n-torus, a contradiction; \medskip\newline
(4) $M^+$ is a torus without an only disc. In order that $M$ be a torus, $M^-$ must be a sphere without a disc. Thus item (4) is similar to item (1); \medskip\newline
(5) $M^+$ is a k-torus, $k\geq 1$ without $n\geq 2$ discs. Then $M^-$ is at least an sphere without $n$ discs and $M$ is not a torus. A contradiction (see item (3)).
\qed \medskip

If $M$ is tangent to an equator $S(v)$ then the points of tangency are
isolated from each other (we remember that $M$ has strictly negative
Gauss-Kronecker curvature, so points of tangency $p$ in $S(v)\cap M$ are nondegenerate saddle points regarding $M$ as a graphic over $S(v)$ near $p$).
So one can carry out a $C^{2,\alpha}$ deformation of $M$
by perturbing it in a neighborhood of its points of tangency in order that $S(v)$ is not tangent
to the deformed torus $\overline M$ anymore. If
the deformation is sufficiently small, we can guarantee that $\overline M$ is embedded, has
strictly negative Gauss-Kronecker curvature and the {\it two-piece property} (see Lemma \ref {ltppopen}).

\begin{remark}\label{rdeform}
Here's a local description of the deformation above:
in a small neighborhood of a point
of tangency $p$ the intersection $S(v)\cap M$ is the union of two  
curves crossing at $p$ at a nonnull angle forming an ``x'' shape (nondegenerate saddle point). 
We denote by $q_1,q_2,q_3, q_4$ the points at the extremities of the ``x''. It may be supposed that the points $q_i$ are in a circle in $S(v)$ with center $p$ such that around the circle $q_1$
lies between $q_2$ and $q_4$, $q_2$ lies between $q_1$ and $q_3$, etc. Thus the ``x'' is formed by smooth arcs connecting $q_1$ to $q_3$ and $q_2$ to $q_4$ crossing at a nonnull angle at the point $p\neq q_i$. We assume the deformation near $p$ restricted
to a small ball $B$ with center $p$ such  that $q_1,\dots,q_4\not \in B$. Since $\overline M$ is not tangent to $S(v)$, locally   $S(v)\cap \overline M$ is the disjoint union of two curves that either connect $q_1$ to $q_2$ and $q_3$
to $q_4$ or connect $q_1$ to $q_4$ and $q_2$ to $q_3$. See in \cite {Ad}, p. 101, a nice picture of this bifurcation.  
\end{remark}

\begin{lemma}\label{lclass2} An equator $S(v)$ can be tangent to $M$ at most at two points.
\end{lemma}

\noindent{\it Proof.} Suppose $S(v)$ tangent to $M$. The assumptions on $M$ imply that
$S(v) \backslash M$ is the union of open regions whose borders are closed curves embedded in $S(v)$.  We remark that these curves are not necessarily smooth, for points of tangency may belong to some of them. 

Let $\mathcal C_1$, $\mathcal C_2$ be the connected components of $\S3 \backslash M$. If  both $\mathcal C_1\cap S(v)$ and $\mathcal C_2\cap S(v)$ have only one connected component, then $S(v)\cap M$ contains an only closed curve embedded in $S(v)$, that is to say, the border of these connected components in $S(v)$. Besides, which is contrary to our assumptions, $S(v)$ is not tangent to $M$ at any point because near points of tangency curves in $S(v)\cap M$ must have the ``x" shape described in Remark \ref{rdeform}. 

Thus we can assume that  $\mathcal C_1\cap S(v)$ has at least two connected components. Firstly, we assume that  $\mathcal C_1\cap S(v)$ has exactly two connected components bordered, respectively, by the closed curves $l_1$ and $l_2$ embedded in $S(v)$ (so both components are discs in $S(v)$).
Observe that 
$l_1\cap l_2$ may have only points where $S(v)$ and $M$ are tangent so it
is a discrete (possibly even empty) set ($S(v)$ is tangent to $M$ at a discrete set of points because $M$ has strictly negative Gauss-Kronecker curvature). In fact, if $l_1\cap l_2$ had a point $q$ where $S(v)$ and $M$ are strictly transversal then a line in $S(v)$ would cross $l_1$ and $l_2$ at $q$  without interchanging connected components of $\S3 \backslash M$, a contradiction since $M$ is supposedly embedded.

Admit that  $l_1\cap l_2$ has three or more points $q_1,q_2,\dots,q_n$ disposed sequentially along $l_1$ (and thus along $l_2$ as can be seen). Thus the union of closed arcs in $l_1$ and $l_2$ with extremities $q_1$ and $q_2$ is a closed embedded curve $\lambda_1$ bounding a disc in $S(v)$. Analogously embedded closed curves $\lambda_2$ and $\lambda_3$ are obtaided from arcs in $l_1$ and $l_2$ connecting $q_2$ to $q_3$ and $q_3$ to $q_4$ (or to $q_1$,
if there are only 3 points in  $l_1\cap l_2$), respectively.

Now we will deform $M$ into an embedded $C^{2,\alpha}$ torus $\overline M$ with the {\it two-piece property} and strictly negative Gauss-Kronecker curvature such that $\overline M$ and $S(v)$ are not tangent (Lemma \ref{ltppopen}).  As seen in Remark \ref{rdeform}, along the deformation the curves can bifurcate at the tangency points so that $\lambda_1$, $\lambda_2$ and $\lambda_3$ are taken to disjoint embedded closed curves $\overline \lambda_1$, $\overline \lambda_2$ and $\overline \lambda_3$, 
respectively. This contradicts  Lemma \ref{lclass1}. 

Now let $\mathcal C_1\cap S(v)$ have two or more connected components  bordered by embedded closed curves $l_1$, $l_2,\dots,l_n$, $n\geq 3$, where as above distinct curves do not cross each other (they border connected components of $\mathcal C_1$) and may only have a discrete set of points in common. In order to finish the proof of this lemma we can either proceed as above or in a more straightforward approach separate the borders of these connected components by a small perturbation of $M$ around the points of tangency between $M$ and $S(v)$, which contradicts  Lemma \ref{lclass1} again.\qed \medskip

Lemmas above are synthetized in the following proposition:

\begin{proposition}\label{pclass}
Let $M$ be a $C^{2,\alpha}$torus embedded in $\S3$ that has  strictly negative
Gauss-Kronecker curvature. Then $M$ has the  {\it two-piece property} if and only if for every equator $S(v)$ the intersection $S(v)\cap M$  
is classified among the four types listed below:
\begin{enumerate}
\item an embedded closed curve,
\item two disjoint  embedded closed curves,
\item two embedded closed curves bounding disjoint discs in $S(v)$ with only one  point in common,
\item two embedded closed curves bounding disjoint discs in $S(v)$ with exactly two points in common.
\end{enumerate}
where a curve in a intersection of type 1 is nullhomotopic in $M$ and the pair of curves in intersections of
types 2, 3 or 4 are homotopic in $M$ to a same generator of $\pi_1(M)$ for any choice of the pair of curves in type 4 intersections (see Remark \ref {rchico}).
\end{proposition}

\begin{remark}\label{rchico} A type 4 intersection may be described as follows: $S(v)\cap M$ is the union of two embedded smooth closed curves $C_1$ and $C_2$ crossing at two points $P$ and $Q$ where these curves are strictly transversal. Then let $l_1$, $l_2$, $l_3$ and $l_4$ be arcs in $S(v)$ with common extremities $P$ and $Q$ such that, seen as pointsets,  $l_1\cup l_3= C_1$, $l_2\cup l_4=C_2$ and $l_i\cap l_j=\{P,Q\}$ when $i\neq j$. Then one has two choices
for the pair of embedded curves bounding disjoint discs in $S(v)$ meeting only at two points: the curves of trace $l_1\cup l_2$ and $l_3\cup l_4$; or the curves of trace $l_1\cup l_4$ and $l_2\cup l_3$.
\end{remark}

\noindent{\it Proof of Prop. \ref{pclass}.} If the intersection $S(v)\cap M$  is of the types above then necessarily $S(v)$
divides $M$ in two connected components. Indeed, the constraints on the shape and homotopy class of the curves in the intersections  completely characterize the topology of $H_+(v)\cap M$ and $H_-(v)\cap M$ up to  a reflection  with respect to $S(v)$ (see appendix \ref{afig}). Thus $M\backslash S(v)$ either is a disjoint union of a disc with a torus without a disc (type 1), or a disjoint union of two tubes (type 2), and so on.

For the reciprocal, that was mostly the  object of Lemmas \ref{lclass1} and \ref{lclass2}, we assume that $M$ has the {\it two-piece property}. What is left to prove concerns the homotopy class of curves of types 3 and 4. But this characterization follows from the proof of Lemma \ref{lclass2}. In fact, it was proved there that if $S(v)$ is tangent to $M$,  then $S(v)\cap M$ is the union of two embedded closed curves $l_1$ and $l_2$ with at least one and at most two points in common bounding disjoint discs in $S(v)$.  Now deform $M$ into an embedded torus $\overline M$ with the {\it two-piece property} and strictly negative Gauss-Kronecker curvature. From Remark \ref{rdeform} and the proof of Lemma \ref{lclass2}, we may assume that $l_1$ and $l_2$
are taken continuosly along the deformation into disjoint closed curves $\bar l_1$ and $\bar l_2$. By Lemma \ref{lclass1}
these curves are homotopic  generators of $\pi_1(\overline M)$, so $l_1$ and $l_2$ are homotopic generators too.  \qed
\medskip

\begin{remark} 
Types 1, 2, 3 or 4 also describe  the possible intersections 
of  an ordinary torus of revolution in $\mathbb R^3$ with a plane to which it is transversal at some point. This fact can be verified with the aid of the  sketches provided in appendix \ref{afig}. 
\end{remark}

We now assume that $S(v)$ is tangent to an embedded torus $M$ with strictly negative Gauss-Kronecker curvature and the {\it two-piece property}. Let $p$ be one of the points where $M$ and $S(v)$ are tangent. 
In a neighborhood of $p$ the intersection $S(v)\cap M$ contains two disjoint curves crossing $p$ 
at a nonnull angle (since $M$ has strictly negative Gauss-Kronecker curvature). This angle is  formed by the asymptotic directions of $M$ at $p$ (see \cite{Sp}). Admit that a geodesic $\gamma$ passes through $p$ without 
pointing to an asymptotic direction. Near $p$, a rotation of $S(v)$ around $\gamma$ splits the curves  in $S(v)\cap M$ into two connected smooth curves with $\gamma$ between them, regardless the sense of the rotation. While one of these curves does not touches $\gamma$ near $p$ the other one must be tangent to $\gamma$
at $p$ (which passes through $p$ depends on the sense of the rotation). 

If $\gamma$ does not pass through the tangency point $p$ then the bifurcation depends on the sense of the rotation and proceeds exactly as described in Remark \ref {rdeform}. Indeed, Remark \ref{rdeform} treats a particular case in which $p$ is a center of the geodesic $\gamma$, i.e., a point that equidists from the geodesic. Thus, in the notation of Remark \ref{rdeform}, if $q_1$, $q_2$, $q_3$ and $q_4$ are the extremities of the ``x" shape in $S(v)\cap M$ with center $p$, then we may assume that a small rotation of $S(v)$ around $\gamma$ to an equator $S(\bar v)$ takes $q_1$, $q_2$, $q_3$ and $q_4$ to points $\bar q_1$, $\bar q_2$, $\bar q_3$ and $\bar q_4$, respectively, in $S(\bar v)\cap M$, such that in a neighborhood of $p$, according to the sense of the rotation, either $\bar q_1$ is connected to $\bar q_2$ and $\bar q_3$ to $\bar q_4$ by disjoint arcs in $S(\bar v)\cap M$,  or $\bar q_1$ is connected to $\bar q_4$ and $\bar q_2$ to $\bar q_3$.

Thus let $\gamma\subset S(v)$ be a geodesic that does not point to asymptotic directions of $M$ at points where $M$ is tangent to $S(v)$. If $S(\bar v)$ is obtained by a small rotation of $S(v)$ around $\gamma$ then $S(\bar v)$ 
is not tangent to $\gamma$, provided the rotation is sufficiently small. Then, as we have seen, near the points of tangency the intersection $S(\bar v)\cap M$ is the union of regular disjoint curves and outside a neighborhood of the points of tangency both $M$ and $S(v)$ (hence also $M$ and $S(\bar v)$) are strictly transversal if the rotation is sufficiently small.

These remarks will be used in the proof of the following result: 

\begin{lemma}\label{ltype} 
There exists an equator $S(v)$ such that $S(v)\cap M$ is of type 2.
\end{lemma}

\noindent
{\it Proof.} 
Firstly, we assume that $M$ is not invariant by the antipodal map. Thus let $p$ be a point in $M$ such that $-p$ does not belong to $M$.
Let $S(\bar v)$ the equator tangent to $M$ at $p$. 

If $S(\bar v)\cap M$ is of type 3, let $D_1$ and $D_2$ be the two disjoint discs in $S(v)$ bounded by closed curves in $S(\bar v)\cap M$ 
whose common point is $p$. 
Then assume that the geodesic $\gamma\subset S(\bar v)$ passes through $p$ without pointing to an asymptotic direction and does not cross $D_1$ and $D_2$. Thus, near $p$, the geodesic $\gamma$ lies 
between $D_1$ and $D_2$. Hence, in a neighborhood of $p$, the intersection $S(v)\cap \gamma$ is the  union of four arcs, namely $l_1$, $l_2$, $L_1$ and
$L_2$, with only common point the extremity $p$ such that $l_1$ and $l_2$ are in the border of $D_1$ and $L_1$, $L_2$ border $D_2$ (see Remark \ref{rdeform}). 
Let $l_3=\partial D_1\backslash (l_1\cup l_2)$ and $L_3=\partial D_2\backslash (L_1\cup L_2)$.

A small rotation of $S(\bar v)$ around $\gamma$ into an equator $S(v)$ takes $(l_1\cup l_2)$ 
and $L_1\cup L_2$ into smooth curves $l_\alpha$, $L_\alpha$ in $S(v)\backslash \gamma$, respectively. As seen above, these
curves are smooth, disjoint and rest in different sides of $\gamma$. As regards the disjoint smooth curves $l_3$ and $L_3$, these curves are taken  by the small rotation into disjoint smooth curves in $S(\bar v)\cap M$, say $l_\beta$ and $L_\beta$. 
We may assume that the extremities of $l_\alpha$ and $L_\alpha$ coincide with the extremities of $l_\beta$ and $L_\beta$, respectively.
Then provided the rotation is sufficiently small, $l_\alpha \cup l_\beta$ and
$L_\alpha \cup L_\beta$ are the trace of disjoint embedded closed curves and $S(v) \cap M$ is of type 2.

If $S(\bar v)\cap M$ is of type 4 with $S(\bar v)$ tangent to $M$ at $p$ and $q$ then $S(\bar v)\backslash M$ contains among its components two disjoint discs $D_1$ and $D_2$ whose boundaries in
$S(\bar v)\cap M$ are two closed curves with common points only $p$ and $q$. 
As above let $\gamma\subset S(\bar v)$ be a geodesic through $p$ that does not points to an asymptotic direction and also separates $D_1$ from $D_2$ near $p$. Since $M$ is not invariant by the antipodal map we may 
assume that $q\neq -p$ so we may also suppose that $\gamma$ does not pass through $q$.We define $l_1$, $l_2$
and $L_1$, $L_2$ in a neighborhood of $p$ as in the preceding argument for type 3 intersections. Analogously we denote by $l_3$, $l_4$, $L_3$, $L_4 \subset S(\bar v)\cap M$ the arcs in a small neighborhood of $q$   with a extremity in $q$  such that
$l_3\cup l_4\subset \partial D_1$ and  $L_3\cup L_4\subset \partial D_2$. The sets $S_1=\partial D_1 \backslash (l_1\cup l_2\cup l_3\cup l_4)$ and $S_2=\partial D_1 \backslash (l_1\cup l_2\cup l_3\cup l_4)$ are the disjoint union of smooth arcs $l_5,$ $l_6\subset \partial D_1$ and $L_5,$ $L_6\subset \partial D_2$, respectively.

After a small rotation of $S(\bar v)$ into $S(v)$ the curves $l_1$, $l_2$, $L_1$, $L_2$ are taken to smooth disjoint curves $l_\alpha$ an $L_\alpha$ of $S(v)\cup M$ as in the previous case. We may choose the sense of the rotation (see
discussion above and Remark \ref{rdeform}) in order  that $l_3$, $l_4$ and  $L_3$, $L_4$ are taken, respectively, to  smooth disjoint arcs $l_\beta$ (obtained from $l_3$, $l_4$) and $L_\beta$ (obtained from $L_3$, $L_4$) as described in Remark \ref{rdeform}. Arcs  $l_5$, $l_6,$ $L_5$ and $L_6,$ are taken to smooth disjoint arcs $l_\mu$, $l_\nu$, $L_\mu$ and $L_\nu,$ respectively. So we can assume that
$l_\alpha\cup l_\beta \cup l_\mu\cup l_\nu$ and $L_\alpha\cup L_\beta\cup L_\mu\cup L_\nu$ are the trace of disjoint closed curves in $S(v)\cap M$ that is thus of type 2.

Finally, let $M$ be invariant by the antipodal map and $S(\bar v)$ be tangent to $M$ at the point $p$. Then $S(\bar v)$ is tangent to $M$ at the point $-p$ too. Thus $S(\bar v)$ is of type 4 and there exists disjoint discs $D_1$and $D_2$ that are components of $S(\bar v)\backslash M$ whose borders have common points $p$ and $-p$. 
Let $\gamma$ be a geodesic of $S(\bar v)$  such that at p $\gamma$ does not point to an asymptotic direction and separates $D_1$ from $D_2$ in a neighborhood of $p$. One may verify that the antipodal map preserves connected components of $\S3\backslash M$. Thus the antipodal invariance of $M$ implies that at
$-p$ the geodesic $\gamma$ does not point to an asymptotic direction and locally separates $D_1$ from $D_2$.
As in the analysis of a type 3 intersection, once a small rotation of $S(\bar v)$ into $S(v)$ is performed, the discs $D_1$ and $D_2$ are set apart into connected components of $S(v)\backslash M$ with disjoint boundaries. So $S(v)\cap M$ is of type 2. \qed

\begin{lemma} \label{ltype2} Let $T$ be the Clifford torus and $S(v)$ an equator tangent to $T$. Let $(v_n)\subset \S3$ be a sequence in $\S3$ such that $v_n\rightarrow v$ and each $S(v_n)$ is not tangent to $T$ 
at any point. Then
\begin{enumerate}
\item each intersection $S(v_n)\cap T$ is of type two;
\item if $\lambda_n$ is a closed curve in $S(v_n)\cap T$ then there exists points $p_n,q_n\in \lambda_n$ such that $k_n(p_n)\rightarrow \infty$ and $k_n(q_n)\rightarrow 0$ when $n\rightarrow \infty$, where $k_n$ is the curvature of $\lambda_n$. 
\end{enumerate}
\end{lemma}

\noindent
{\it Proof.}\smallskip

\noindent
(1) $\S3\backslash T$ is the union of two connected components that are invariant by the antipodal map. Now assume that $S(v_n)\cap T$ is of type 1. The connected components $S_1,$ $S_2$ of $S(v_n)\backslash T$ are in different connected components of $\S3 \backslash T$. Moreover, the antipodal map takes $S_1$ to $S_2$ and vice versa. A contradiction.\smallskip

\noindent (2) In the proof of Lemma \ref{ltype} it was argued that  the invariance of $T$ by the antipodal map 
implies that $S(v)\cap T$ is of type 4. Besides, $S(v)\cap T$ has all the compatible symmetries. A curve in $S(v)\cap T$ separating two regions must bend equally into the two regions, so its curvature must be equal to zero. The intersection   $S(v)\cap T$ is then the union of two geodesics  crossing at the antipodal points $p$, $-p$ at an angle equal to $\pi/2$. By taking a subsequence, if necessary, we may assume that the sequence $(\lambda_n)$ converges to arcs of geodesics in $S(v)\cap T$ bounding a component of $S(v)\backslash T$. One may choose points $q_n\in \lambda_n$ away from $p$, $-p$ at least a sufficiently small $\epsilon>0$. Thus $k_n(q_n)\rightarrow 0$. Assume that the curvature of $\lambda_n$ remains bounded when $n\rightarrow \infty$. Thus  a subsequence of $(\lambda_n)$ converges to a curve $\lambda$ that at least belongs to $C^{1,\alpha}$ for $0\leq\alpha<1$ (since $C^2(S^1)$ is compact in $C^{1,\alpha}(S^1)$, see Theorem 1.31 in \cite{Am}). A contradiction, because $(\lambda_n)$ converges to arcs of geodesic that meet at an angle equal to $\pi/2$. So there exists the sequence $(p_n)$ defined above. \qed

\begin{lemma}\label{ltype3} Let $M$ be a minimal torus embedded in $\S3$ that is not congruent to the Clifford torus $T$. Then there exists an equator $S(v_0)$ such that $S(v_0)\cap M$ is of type 2 and contains a closed curve that is not congruent to any closed curve in a intersection $S(v)\cap T$ of type 2.
\end{lemma}

\begin{remark}\label{rtype3} Lemma \ref{ltype3} seems quite straightforward since $M$ and $T$ are noncongruent real analytic tori (see Lemma 1.1 in \cite{L1}).  Its proof below would be simpler had we used a fact communicated to us by H. Rosenberg, 
namely, that whenever a torus minimally embedded in $\S3$ contains a geodesic it is congruent to the Clifford torus.  
\end{remark}

\noindent
{\it Proof of Lemma \ref{ltype3}.} Firstly, we recall that the torus $M$ has strictly negative Gauss-Kronecker curvature and the {\it two-piece property} (by Lemma 1.4 in \cite{L1} and Theorem 2 in \cite{Rs}; see also Int.). We also observe that there exists an equator $S(w)$ such that $S(w)\cap M$ is of type 2 (Lemma \ref{ltype}). Suppose that whenever $S(w)\cap M$ is of type 2 then both closed curves in $S(w)\cap M$ are  congruent to closed curves in intersections $S(\bar w_1) \cap T$ and $S(\bar w_2)\cap T$, respectively. 

We may choose $v_n\in \S3$ such that $S(v_n)\cap M$ is of type 2 for every $n$. Moreover, we can also assume that  $v_n \rightarrow v$, and that $S(v)$ is tangent to $M$ (just take $v_0\in \S3$ such that $S(v_0)\cap M$ is of type 2 and move $S(v_0)$ until it becomes tangent to $M$).  

Thus $S(v_n)$ converges to $S(v)$ and  curves in $S(v)\cap M$ must be the limit (in the sense discussed in item 2 of the proof of Lemma \ref{ltype2}) of closed curves congruent to curves in intersections $S(w)\cap T$. Besides curves in $S(v)\cap M$ are not in $C^2$ because of their singularities at the tangency  points. Thus {\it $S(v)\cap M$ must be congruent to an intersection $S(\bar v)\cap T$ of type 4 such that $S(v)\cap M$ is the intersection of two geodesics of $\S3$ crossing at an angle equal to $\pi/2$} (see proof of Lemma \ref{ltype2}).

We consider then  the set $S$ of points $p$ where $M$ is tangent to equators $S(v_p)$ such that $S(v_p)\cap M$ is the union of two geodesics crossing at an angle equal to $\pi/2$. This set is readily seen to be closed and not empty, by Lemma \ref{ltype2} and the remarks above. Let $p \in S$ and $V$ be a small neighborhood  of $p$ in $M$. Let $S(v_q)$ be an equator tangent to $M$ at some point $q$ in $V$. We have the following situations to analyze:

\noindent
(1) If $S(v_q)\cap M$ is of type 3 then we proceed as in the proof of Lemma \ref{ltype}. We rotate $S(v_q)$ around a geodesic $\gamma\subset S(v_q)$ through $q$ that lets the embedded closed curves in $S(v_q)\cap M$ in the closures of different connected components of $S(v_q)\backslash \gamma$. As was seen above the rotated equators have intersections of type 2 with $M$ whenever the rotation is sufficiently small. Our hypothesis thus implies that $S(v_q)\cap M$ is the union of two geodesics meeting  at two points at an angle equal to $\pi/2$, a contradiction.

\noindent
(2) Now suppose that $S(v_q)\cap M$ is of type 4. If $S(v_p)$ is the equator tangent to $M$ in $p$
then there exists a geodesic $\gamma$ in $S(v_p)$  such that the closure of each component of $S(v_p)\backslash \gamma$ contains an embedded closed curve. Hence when $V$ is sufficiently small and $q\in V$ the same happens to any intersection
$S(v_q)\cap M$ of type 4. Thus we rotate the equator $S(v_q)$ around a geodesic in  corresponding to $\gamma$ above in order to obtain type 2 intersections as in the previous case. Then we proceed as in item 1 above.

So $S$ is open and closed in $M$ besides being not empty. Hence $S=M$ and the asymptotic curves of $M$ (i.e., curves whose normal curvature is equal to zero, see \cite{Sp}, vol. 3) are geodesics of $\S3$. 
Then $M$ is an embedded flat minimal torus (quick proof: if $p\in M$, consider a sequence $R_n$ of rectangles with vertex $p$ whose sides are arcs of asymptotic curves; then assume that the length of the sides of $R_n$ tends to zero when $n\to\infty$ and apply the Gauss-Bonnet formula). According to a well known result due in \cite{cdck} and \cite{L2}, the minimal torus $M$ must be the Clifford torus, a contradiction.\qed  \medskip

The next result is quite straightforward as well and may also be obtained from the characterization of the Clifford torus in Remark \ref{rtype3}.

\begin{lemma}\label{ltype4} Let $M$ be a torus minimally embedded in $\S3$ and assume that $S(v)\cap M$ is of type 2. Then none of the closed curves in $S(v)\cap M$ is a geodesic of $\S3$.
\end{lemma}

\noindent
{\it Proof.} Let us admit that $\lambda_1\subset S(v)\cap M$  is a geodesic of $\S3$. We perform a rotation of $S(v)$ around $\lambda_1$ until an equator $S(w)$ obtained from the rotation of $S(v)$ becomes tangent to $M$
at a point $p$ (such equator does exist because $\lambda_1\subset M$). If $p \not \in\lambda_1$ then $S(w)\cap M$ is not of types 1, 2, 3 or 4 and $M$ does not have the  
{\it two-piece property} (see Prop. \ref{pclass}), which contradicts Ros' Theorem (i.e., Theorem 1.2 in \cite{Rs}). 
Thus $p\in \lambda_1$. 

Now we observe that $S(v)\cap M$ is of type 2. Hence $S(v) \backslash \lambda_1$ is the union 
of two discs $D_1$ and $D_2$ and we may assume from the definition of type 2 intersections  that $D_1\cap M=\emptyset$ 
and that $D_2\cap M$ is a closed curve. Along the rotation, until the first point of tangency, the intersections of the equators with $M$ remain of type 2. Hence the disc rotated from $D_2$
contains a closed curve and the disc rotated from $D_1$ has empty  intersection with $M$. Let the discs $D_1,\ D_2$ be taken by the rotation above to discs $D_1^*,\ D_2^*\subset S(w)$, respectively, which are the components of $S(w)\backslash \lambda_1$. Thus,
$D_1^*\cap M=\emptyset$ (or $M$ does not have strictly negative Gauss-Kronecker curvature) and, consequently,  the point of tangency $p$ cannot be a nondegenerate saddle point (a curve through $p$ in $S(w)\cap M$  must be in the closure of $D_2^*$). This contradicts the already cited fact (for example, in the proof of the preceding lemma) that $M$ has strictly negative Gauss-Kronecker curvature. \qed

\section{Sequences of minimal tori converging to the Clifford torus}\label{scon}

Let the open annulus $A_2=\{x\in \mathbb R^4:\hspace{.03in} 1/2<|x|<2\}$ and $F:A_2\to A_2$ be a $C^{2,\alpha}$ map.  We let $D^j=\partial /\partial x_j$, $D^{jk}=\partial^2/\partial x_j\partial x_k$ and set
\begin{eqnarray}\label{eqnh1}
[F]_{\alpha}&=&\sup_{{x,y\in A_2}\atop{x\neq y}}\frac{|F(x)-F(y)|}{|x-y|^\alpha},\\ \nonumber
||F||_{C^{2,\alpha}}&=& \sup_{A_2}|F(x)|+ \max_{1\leq j \leq 4} \sup_{A_2}\left|{D^j F(x)}\right|+\max_{1\leq j,k \leq 4} \sup_{A_2}\left|D^{jk} F(x)\right|\\
&&+\max_{1\leq j,k \leq 4}[D^{jk} F]_{\alpha},
\end{eqnarray}
a seminorm and the Holder $C^{2,\alpha}$ norm in $A_2$, respectively, where $|\hspace{.06in}|$ is the euclidean norm in $\mathbb R^4$.

\begin{remark}\label{remb} If $F$ can be extended in $C^2$ to a map $G$ defined in  the closure $A_2^\circ$ then it is easily seen that $||F||_{C^{2,\alpha}}$ and  $||G||_{C^{2,\alpha}(A_2^\circ)}$ coincide.
\end{remark}

The purpose of the following lemma is to serve as a compacity criterion for what follows. Let $X,\ Y,$ where $X\subset Y$, be topological spaces. We observe that only in Lemma \ref{lemb} below the term ``embedding'' refers to the identity map $I:X\to Y$ when $I$ is continuous. If there also holds that $I(X)$ is precompact in $Y$ then we say that $X$ is compactly embedded in $Y$. 

Any $C^{2,\alpha}$ diffeomorphism $\xi$ of $\S3$ is canonically 
extended to a $C^{2,\alpha}$ diffeomorphism $X$ of $A_2$: if $v\in \S3$ and $1/2< r<2$ let $X(rv)=r\xi(v)$. Thus

\begin{lemma}\label{lemb} Let $\mathcal{X}^\alpha$ be the set of  $C^{2,\alpha}$ diffeomorphisms of $\S3$ canonically extended to $A_2$. If  $0<\nu<\lambda\leq 1$  then the embedding $\mathcal{X}^\lambda\to \mathcal{X}^\nu$ exists and is compact. 
\end{lemma} 

\noindent
{\it Proof.} Lemma \ref{lemb} follows directly from Remark \ref{remb} and Theorem 1.31 in \cite{Am}. \qed  \medskip

We let now $I:A_2\to A_2$ be the identity map in $A_2$ and define the following functional on the diffeomorphisms of $A_2$:
\begin{equation}
\tau^\alpha(X)= ||X-I||_{C^{2,\alpha}}+||X^{-1}-I||_{C^{2,\alpha}},
\end{equation}
where $0<\alpha\leq 1$.

\begin{lemma}\label{lbas} If $0<\alpha\leq 1$ let $(X_k)$ be a sequence of $C^{2,\alpha}$ diffeomorphisms of $\S3$ canonically extended to $A_2$ as above. Assume that $\tau^\alpha(X_k)\leq C<\infty\ \forall\hspace{.02in}k\in \mathbb N$.  Let  $0<\beta<\alpha$. Then a subsequence of $(X_k)$ converges  in
$C^{2,\beta}$ to a diffeomorphism $X$ of $\S3$ canonically extended to $A_2$. Moreover,  besides being a map in $C^{2,\beta}$, $X$ is also a $C^{2,\alpha}$ diffeomorphism and ${\tau^{\alpha}(X)}\leq C$. Lastly,  if $X_k \to X$ in $C^{2,\alpha}$, then $X$ is a $C^{2,\alpha}$ diffeomorphism of $\S3$ canonically extended to $A_2$, and there exists $\lim_{k\to\infty}\tau^\alpha(X_k)=\tau^\alpha(X)$.
\end{lemma}

\noindent
{\it Proof.} Since the norms $||X_k-I||_{C^{2,\alpha}}$ are bounded for all $k$, the sequence $(X_k)$ is bounded in $C^{2,\alpha}$. By Lemma \ref{lemb} there exists a subsequence, denoted $(X_{k_l})$, that converges in $C^{2,\beta}$ to a map $X:A_2\to A_2$. Now we consider  the sequence $(X_{k_l}^{-1})$ obtained from the subsequence above. From the fact that the norms $||X_{k_l}^{-1}-I||_{C^{2,\alpha}}$ are also bounded for all $k$ and Lemma \ref{lemb} we have that a subsequence of $(X_{k_l}^{-1})$, denoted $(X_{k_l}^{-1})$ as well, converges in $C^{2,\beta}$ to a $C^{2,\beta}$ map $X^*:A_2\to A_2$. A continuity argument implies that exists  $X^{-1}=X^*$. Hence $X$ is a  $C^{2,\beta}$ diffeomorphism of $\S3$ canonically extended to $A_2$.

In fact, $X$ is a $C^{2,\alpha}$ diffeomorphism.  By a routine argument, let $1\leq i,j \leq 4$: since $\lim_{l\to\infty}D^{ij} X_{k_l}(x)=D^{ij} X(x)$ then
for every pair $x,y\in A_2,\ x\neq y,$ there holds that 
$$\lim\inf_{l\rightarrow\infty} [D^{ij}(X_{k_l}-I)]_{\alpha}\geq
\lim_{l\rightarrow\infty}\frac{|D^{ij} (X_{k_l}-I)(x)-D^{ij}(X_{k_l}-I)(y)|}{|x-y|^\alpha}=\hspace{0.5in}$$
$$\hspace{1.8in}\frac{|D^{ij} (X-I)(x)-D^{ij}(X-I)(y)|}{|x-y|^\alpha}.$$ 
Thus both $[D^{ij}(X-I)]_\alpha$, $||X||_{C^{2,\alpha}}<\infty$. In fact,  $\lim\inf_{l\rightarrow\infty} [D^{ij}(X_{k_l}-I)]_{\alpha}\geq [D^{ij}(X-I)]_\alpha$.
Besides, since it was assumed  that $X_{k_l}\rightarrow X$ in $C^{2,\beta}$ for some $0<\beta< \alpha$, the other terms in the expression of $||X_{k_l}-I||_{C^{2,\alpha}}$ converge to their counterparts in the expression of $||X-I||_{C^{2,\alpha}}$. From these remarks there follows that
$\lim\inf_{l\rightarrow \infty}||X_{k_l}-I||_{C^{2,\alpha}}\geq ||X-I||_{C^{2,\alpha}}$. Analogously we have the corresponding inequality  for the inverse map: $\lim\inf_{l\rightarrow \infty}||X_{k_l}^{-1}-I||_{C^{2,\alpha}}\geq||X^{-1}-I||_{C^{2,\alpha}}$. Thus $C\geq\lim\inf_{l\to\infty}\tau^\alpha(X_{k_l})\geq\tau^\alpha(X)$. 

If $X_k \to X$ in $C^{2,\alpha}$ then $(X_k-I)$ and $X_k^{-1}-I$ converge to $X-I$ and $X^{-1}-I$ in $C^{2,\alpha}$, respectively, since $\tau^\alpha(X_k)$ is bounded, and the Lemma follows.\qed

\medskip

The Clifford torus $T$ (see  Def. \ref{dctorus2}) is reticulated by the longitudes and meridians below:
\begin{eqnarray}
\phi^k_n(\theta)=\frac{\sqrt{2}}{2}(\cos(\frac{\pi k}{n}),\sin(\frac{\pi k}{n}),\cos(\theta),\sin(\theta))\\
\psi^k_n(\theta)=\frac{\sqrt{2}}{2}(\cos(\theta),\sin(\theta),\cos(\frac{\pi k}{n}),\sin(\frac{\pi k}{n}))
\end{eqnarray}
where $\theta \in [0,2\pi)$ and $0\leq k<2n$.  The components of $T$ minus the longitudes and meridians above are the $4n^2$ squares in $T$ with vertices the points
\begin{equation}
p^{jk}_n=\frac{\sqrt{2}}{2}(\cos(\frac{\pi j}{n}),\sin(\frac{\pi j}{n}),\cos(\frac{\pi k}{n}),\sin(\frac{\pi k}{n})).
\end{equation}
Now let $v_0=(0,1,0,0)$.
Thus $S(v_0)\cap T$ is of type 2 and is the disjoint union of the lines of curvature $\phi^0_n$ and $\phi_n^n$.

Hereafter we assume the curves in $S(v_0)$ parameterized by arclength. So 
let $\lambda$ be a $C^{2,1}$ closed curve embedded in the equator $S(v_0)$ and $\Omega_n^\lambda$ be the set of  $C^{2,1}$ diffeomorphisms $X$ of $\S3$ canonically extended to $A_2$ that let $S(v_0)$ invariant such that both $X(\phi_n^0)$ is congruent to $\lambda$ (seen as a pointset) and the image $X(T)$
of the Clifford torus has mean curvature equal to zero at the points $X(p^{jk}_n)$ (abusing of the notation, here $T$, $\lambda$, $\phi_n^0$, etc., when arguments of maps, must be regarded as pointsets). 

\begin{remark}\label{rconst} We observe that $\Omega_n^\lambda$ is not empty. In fact, it is seen that there exists 
 a diffeomorphism $\xi$ of $S(v_0)$ taking $\phi^0_n$ to $\lambda$. Since the set of the $p_n^{jk}$ is discrete, we may assume, by perturbing the diffeomorphism 
$\xi$, if necessary, that $\xi(T)$ has mean curvature zero at the points $\xi(p_n^{jk})$. This diffeomorphism is then extended canonically to a diffeomorphism $X$ of $A_2$ as described above. 
\end{remark}

We now define some functions on the set of closed $C^{2,1}$ curves $\lambda$ embedded in $S(v_0)$. Firstly, we set
\begin{equation}
\delta_n(\lambda)=\inf\{\hspace{.02in}\tau^1(X)\hspace{.02in}:\hspace{.02in}X\in \Omega_n^\lambda\}.
\end{equation}
Thus $\delta_n(\lambda)=0$ if and only if $\lambda$ is congruent to $\phi_n^0$.

If  $0<\alpha<1$ and $M$ is an unknotted $C^{2,\alpha}$ torus (see \cite{L3}) embedded in $\S3$ let $\Lambda_{\alpha}^M$ be the set of $C^{2,\alpha}$ diffeomorphisms $Y$ of $\S3$ canonically extended to $A_2$ such that $Y(T)$ is congruent to $M$. 
Let now $\Omega_{n,\alpha}^{\lambda,\mathcal{K}}$ be the closure in $C^{2,\alpha}$ of the set of diffeomorphisms $X\in \Omega_n^\lambda$ such that $\tau^1(X)<\mathcal{K}$. We remark that $\Omega_{n,\alpha}^{\lambda,\mathcal{K}}\subset \Omega_n^\lambda$ by Lemma \ref{lbas}.
Let $0<\alpha<1$ and $\mathcal{K}> \delta_n(\lambda)$, and set
\begin{equation}
\xi_{n,\alpha}^\mathcal{K}(\lambda)=\inf\{\tau^\alpha(Y)\hspace{.02in}:\hspace{.02in} Y\in \Lambda^{X(T)}_\alpha
\hbox{ for some }\hspace{.02in}X\in\Omega_{n,\alpha}^{\lambda,\mathcal{K}}\}. 
\end{equation}

Given an arclength parameterized arc $\mu$ with trace in $S(v_0)$ we obtain a closed arclength parameterized curve $\lambda$ in $S(v_0)$ by joining the extremities of $\mu$ by an arclength parameterized arc $\nu$ in $S(v_0)$ and letting $\lambda$ 
be the arclength parameterized curve whose trace is the union of the traces of $\mu$ and $\nu$. We may suppose that $\lambda$ is in $C^{2,1}$ by assuming that both $\mu,\nu \in C^{2,1}$ and that $\mu, \nu$ join regularly (i.e., $C^{2,1}$) at
 their extremities. We denote this construction by the product $\lambda=\mu\nu$, where it is transparent that we identify distinct arclength parameterizations of $\lambda$.

Let $\lambda$ be a closed $C^3$ curve embedded in $S(v_0)$ parameterized by arclength and 
$\mathcal{K}\geq\delta_n(\lambda)+1$. 
Since $\lambda$ is closed, we may consider it parameterized for all $s$ by letting
$\lambda(s)=\lambda(s\mod 2L)$, where $2L$ is its length. 
We 
 set $\lambda_t=\lambda|_{[-(1-t)L, (1-t)L]}$, 
$0< t< 1$, for the restriction of $\lambda$ to the interval $[-(1-t)L, (1-t)L]$. Let
 $\mathcal{S}_{n,\alpha}^t$, $0<\alpha,\ t<1$, be the set of $C^{2,1}$ arcs $\mu$ with the same extremities as $\lambda_t$
such that the product $\lambda_t\mu$ is a closed $C^{2,1}$ curve embedded in $S(v_0)$ and 
$\xi_{n,\alpha}^{\mathcal{K}}(\lambda_t\mu)$ is defined. This set is not empty because it contains the arc of $\lambda$ whose trace is
$\lambda\backslash \lambda_t$. We set

\begin{equation}\label{eqphi}
\Phi_{n,\alpha}^{\lambda}(t)=\inf_{\mu\in \mathcal{S}_{n,\alpha}^t}\xi_{n,\alpha}^{\mathcal{K}}(\lambda_t\mu)
\end{equation}
if $t\in (0,1)$ and let $\Phi_{n,\alpha}^{\lambda}(0)=\xi_{n,\alpha}^{\mathcal{K}}(\lambda)$.

\begin{lemma}\label{lcon1} The function $\Phi_{n,\alpha}^{\lambda}$ is continuous and nonincreasing.
\end{lemma}

\noindent
{\it Proof.} Let $t_2>t_1$ and $\mu \in \mathcal{S}_{n,\alpha}^{t_1}$. Thus $(\lambda_{t_1}\mu)\backslash\lambda_{t_2}\in \mathcal{S}_{n,\alpha}^{t_2}$ and, consequently, $\Phi_{n,\alpha}^{\lambda}$ is nonincreasing. Let  $(\mu_k^t)\subset\mathcal{S}_{n,\alpha}^t$ be a minimizing sequence for $\Phi_{n,\alpha}^{\lambda}(t)$, i.e.,
$\lim_{k\to \infty}\xi_{n,\alpha}^{\mathcal{K}} (\lambda_t\mu_k^t) =\Phi_{n,\alpha}^{\lambda}(t)$.  
Taking into account that $\Phi_{n,\alpha}^{\lambda}$ is nonincreasing we admit that
$\lim_{t\stackrel{-}{\to}t_0}\Phi_{n,\alpha}^{\lambda}(t)>\Phi_{n,\alpha}^{\lambda}(t_0)$ and let
$2\epsilon=\lim_{t\stackrel{-}{\to}t_0}\Phi_{n,\alpha}^{\lambda}(t)-\Phi_{n,\alpha}^{\lambda}(t_0)$.
If $k$ is sufficiently large then $\xi_{n,\alpha}^{\mathcal{K}} (\lambda_{t_0}\mu_k^{t_0})-\Phi_{n,\alpha}^{\lambda}(t_0)<\epsilon$. Thus there exists 
$(X,Y)\in \Omega_{n,\alpha}^{\lambda_{t_0}\mu_k^{t_0},\mathcal{K}}\times \Lambda^{X(T)}_{\alpha}$ such that 
$\tau^\alpha(Y)-\Phi_{n,\alpha}^{\lambda}(t_0)<\epsilon$. Since $\Omega_{n,\alpha}^{\lambda_{t_0}\mu_k^{t_0},\mathcal{K}}$ is the closure in $C^{2,\alpha}$
of $C^{2,1}$ maps $Z\in \Omega_{n}^{\lambda_{t_0}\mu_k^{t_0}}$ such that $\tau^1(Z)<\mathcal{K}$, then 
we may assume that $\mathcal{K}>\tau^1(X)$. Now we recall that $\lambda$ is a $C^3$ curve:
thus $\lambda_{t_0-\delta}\to\lambda_{t_0}$ in $C^{2,1}$ (see Remark \ref{rtopos}) and we may admit for each $\delta$ 
the existence of  $\mu_\delta\in \mathcal{S}_{n,\alpha}^{t_0-\delta}$ such that 
$\lambda_{t_0-\delta}\mu_\delta\to \lambda_{t_0}\mu_k^{t_0}$ in $C^{2,1}$ when $\delta\to 0$ and thus also in $C^{2,\alpha}$. A continuity argument
implies that if $\delta$ is sufficiently small then there are diffeomorphisms $(X_\delta,Y_\delta)\in \Omega_{n,\alpha}^{\lambda_{t_0-\delta}\mu_\delta,\mathcal{K}}\times \Lambda^{X_\delta(T)}_{\alpha}$ 
(since we may choose $X_\delta$ in order that $\tau^1(X_\delta)<\mathcal{K}$) near $(X,Y)$ in $C^{2,1}\times C^{2,\alpha}$
such that $\tau^\alpha(Y_\delta)-\Phi_{n,\alpha}^{\lambda}(t_0)<\epsilon$, hence $\lim_{t\stackrel{-}{\to}t_0}\Phi_{n,\alpha}^{\lambda}(t)-\tau^\alpha(Y_\delta)>\epsilon$. 
On the other hand, $\tau^\alpha(Y_\delta)\geq \xi_{n,\alpha}^{\mathcal{K}} (\lambda_{t_0-\delta}\mu_\delta)\geq
\Phi_{n,\alpha}^{\lambda}(t_0-\delta)\geq\lim_{t\stackrel{-}{\to}t_0}\Phi_{n,\alpha}^\lambda(t)$, a contradiction.
Thus $\lim_{t\stackrel{-}{\to}t_0}\Phi_{n,\alpha}^{\lambda}(t)=\Phi_{n,\alpha}^{\lambda}(t_0)$. 
{\it Mutatis mutandis}, an slight modification of the argument above  proves that $\lim_{t\stackrel{+}{\to}t_0}\Phi_{n,\alpha}^{\lambda}(t)=\Phi_{n,\alpha}^{\lambda}(t_0)$ even when $t_0=0$. \qed
\medskip

Let $M$ be a minimal torus embedded in $\S3$ that is not congruent to the Clifford torus $T$. We may assume that 
$S(v_0)\cap M$ is of type 2 (by Lemma \ref{ltype} associated with a congruence of $\S3$). Since $S(v_0)\cap M$ is real analytic so are the curves in $S(v_0)\cap M$. Thus, by Lemma \ref{ltype3}, we may assume that
there exists a closed real analytic curve $\lambda^M\subset S(v_0)\cap M$ such 
that $\lambda^M$ is not congruent to any closed curve in a intersection $S(v)\cap T$. 
Moreover, we also know that $\lambda^M$ is not a geodesic of $\S3$ (see Lemma \ref{ltype4}). 

Let $X_M$ be a $C^{2,1}$ diffeomorphism of $\S3$ extended to $A_2$ that lets $S(v_0)$ invariant such that $X_M(\phi^0_n)$ is congruent to $\lambda^M$ and $X_M(T)$ is congruent to $M$. Such a diffeomorphism does exist, since $M$ is unknotted (cf. \cite{L3})  and $\lambda^M$ is a $C^{2,1}$ embedded canonical generator of $\pi_1(M)$. Thus $X_M\in \Omega_n^{\lambda^M}$ 
and $\delta_n(\lambda^M)\leq\tau^1(X_M)$ for all $n\in \mathbb N$.

We choose an arclength parameterization of $\lambda^M$ in order that the curvature of $\lambda^M$ assumes its maximum  $k_{\max}$ (which is necessarily positive since $\lambda^M$ is not a geodesic) at $\lambda^M(0)$. As regards the Clifford torus $T$, there 
exists $v\in\S3$ such that a real analytic embedded closed curve $\mu^T\subset S(v)\cap T$ has a critical point of curvature 
$k_{\max}$ (Lemma \ref{ltype2}) at $\mu^T(0)$ for a suitably chosen arclength parameterization of $\mu^T$. Thus there also exists a diffeomorphism $X_{T}$ such that $X_{T} \in 
\Omega_n^{\mu^T}$ for all $n$ since we may choose $X_{T}(T)$ 
congruent to $T$. Hence $\xi_{n,\alpha}^{\mathcal{K}}(\mu^T)=0$ $\forall\ n\in\mathbb N$ and $\mathcal{K}>0$.

For lemmas \ref{lcon2} and \ref{lfinal} we let $0<\alpha<1$ as usual and set $$\mathcal{K}=\max(\tau^1(X_M),\tau^1(X_T))+1,$$
thus $\mathcal K>\delta_n(\lambda^M)$ as required.

\begin{lemma}\label{lcon2} Let $M$ be a minimal torus embedded in $\S3$ that is not congruent to the Clifford torus  and $\lambda^M$ be the closed curve defined above. Then  there exists a constant $C>0$ 
such that $\Phi_{n,\alpha}^{\lambda^M}(0)\geq C$ $\forall\ n$. Besides, regardless $n$,   $\lim_{t\to 1}\Phi_{n,\alpha}^{\lambda^M}(t)=0$. 
\end{lemma}

{\proof\ } 
Suppose that $\inf_{n\in \mathbb N}\Phi_{n,\alpha}^{\lambda^M}(0)=0$. Then either (1) there exists a subsequence
$(n_k)$ such that $\lim_{k\to \infty}\Phi_{n_k,\alpha}^{\lambda^M}(0)=0$  or (2) $\Phi_{n_0,\alpha}^{\lambda^M}(0)=0$ for some $n_0\in \mathbb N$.
If only (2) happens just consider $n>n_0$ in the hypothesis of the lemma, though the reasoning below also applies to this case.

Thus we may admit that (1) occurs and  then take  maps $(X_{n_k},Y_{n_k})\in \Omega_{n_k,\alpha}^{\lambda^M,\mathcal{K}}\times \Lambda^{X_{n_k}}_{\alpha}$ such that
 $0\leq\tau^\alpha(Y_{n_k})-\Phi_{n_k,\alpha}^{\lambda^M}(0)<1/k$. Thus by Lemma \ref{lbas} there holds  that for 
$0<\beta<\alpha$ subsequences of  $(X_{n_k})$ and $(Y_{n_k})$ converge in $C^{2,\alpha}$ and $C^{2,\beta}$ to diffeomorphisms $X$, $Y\in C^{2,1}\times C^{2,\alpha}$, respectively. By Lemma \ref{lbas}, $0=\lim_{k\to\infty}\tau^\alpha(Y_{n_k})\geq\tau^\alpha(Y)\geq 0$. Thus $Y$ is the identity map.
Since $X_{n_k}(T)$ is congruent to $Y_{n_k}(T)$, there holds that $X(T)$ is congruent to the Clifford torus.
Moreover, $S(v_0)\cap X(T)$ contains a closed curve congruent to $\lambda^M$, which contradicts the definition of $\lambda^M$.

We will now prove that $\lim_{t\to 1}\Phi_{n,\alpha}^{\lambda^M}(t)=0$. 
Let $\lambda^M_t$, $\mu^T_t$ be arcs of the curves $\lambda^M$, $\mu^T$, respectively, obtained from $\lambda^M$, $\mu^T$ as described in the proof of Lemma \ref{lcon1}.
In order to prove that  $\lim_{t\to 1}\Phi_{n,\alpha}^{\lambda^M}(t)=0\ \forall\ n$ we verify the existence of arcs $\nu_t$
such that the product $\lambda^M_t\nu_t\to\mu^T$ in $C^{2,1}$ when $t\to 1$. Indeed, both arcs $\lambda_t^M$, 
$\mu^T_t$ are $C^3$ curves with curvature equal to $k_{\max}$ at $\lambda^M_t(0)$, $\mu^T_t(0)$, respectively, which are also critical points of the curvature function along their respective curves. Thus, up to congruences, the distance from $\lambda^M_t$
to $\mu_t^T$ in the $C^{2,1}$ space of those arcs of curves (that is known to be metrizable, see \cite{Hi}, chapter 1)  tends to zero when $t\to 1$. So for $t$ sufficiently close to 1 we can join  the extremities of $\lambda^M_t$ with an arc $\nu_t$ near   $\mu^T\backslash \mu^T_t$ in $C^{2,1}$  such that $\nu_t\to \mu^T\backslash \mu^T(0)$ in $C^{2,1}$ when $t\to 1$, and the convergence $\lambda^M_t\nu_t\to\mu^T$ in $C^{2,1}$ follows. Hence, there exists diffeomorphisms $X_t$ converging to $X_T$ in $C^{2,1}$ when $t\to 1$  such that $X_t(T)\cap S(v_0)$ contains a closed curve  congruent
to $\lambda_t^M\nu_t$ (just perturb $X_T$ in a neighborhood of $\phi_n^0$). Consequently, $\lim_{t\to 1}\tau^1(X_t)= \tau^1(X_T)<\mathcal{K}$.
 Now we observe  that $\xi_{n,\alpha}^{\mathcal{K}}(\mu^T)=0$ (since $X_T(T)$ is congruent to the Clifford torus  
so the identity map $I\in \Lambda^{X(T)}_{\alpha}$; besides, $\tau^1(X_T)<\mathcal{K}$ from the definition of $\mathcal{K}$). Hence  we can assume the existence of a $C^{2,\alpha}$ diffeomorphism $Y_t$ of $\S3$ canonically extended to $A_2$ associated with $X_t$ that converges to the identity map $I$ in $C^{2,\alpha}$ when $t \to 1$ such that $Y_t(T)$ is congruent to $X_t(T)$ (e.g., one can take $Y_t=X_t\circ X_T^{-1}$). Thus from lemmas \ref{lbas}  and \ref{lcon1} there holds  that $$0\leq \lim_{t \to 1} \Phi_{n,\alpha}^{\lambda^M}(t)\leq \lim_{t\to 1}\xi_{n,\alpha}^{\mathcal{K}}(\lambda_t^M\nu_t)\leq  \lim_{t\to 1}\tau^\alpha(Y_t)=\tau^\alpha(I)=0,$$
and the lemma is proved. \qed

We  then extend continuously $\Phi_{n,\alpha}^{\lambda^M}$ to the interval $[0,1]$ by letting $\Phi_{n,\alpha}^{\lambda^M}(1)=0$. 

\begin{lemma}\label{lfinal} Let $M$ be a minimal torus embedded in $\S3$ that is not congruent to the Clifford torus $T$ and $C$ be the constant provided by Lemma \ref{lcon2}. 
Then for every $\rho \in [0,C]$ there exists a torus $M_\rho$ minimally
embedded in $\S3$ that is congruent to $T$ only when $\rho=0$, and a $C^{2,\alpha}$ diffeomorphism $Y_\rho$ of $\S3$ canonically extended to $A_2$ satisfying $\tau^\alpha(Y_\rho)=\rho$ such that  $Y_\rho(T)$ is congruent to $M_\rho$. 
\end{lemma}

\noindent
{\it Proof.} Let $\lambda^M$ be 
the closed curve embedded in $\S3$ defined above.  By Lemma \ref{lcon2} and the definition of the function $\Phi_{n,\alpha}^{\lambda^M}$   there holds that $\Phi_{n,\alpha}^{\lambda^M}(t_n^\rho)=\rho$
for some $t_n^\rho\in[0,1]$. So for every 
$n \in \mathbb N$ and $\rho\in[0,C]$ there exists  diffeomorphisms $(X_{n,\rho},Y_{n,\rho})\in  \Omega_{n,\alpha}^{C_n^\rho,\mathcal{K}}\times\Lambda^{X_{n,\rho}(T)}_{\alpha}$, 
where $C_n^\rho=\lambda^M_{t_n^\rho}\mu_{t_n^\rho}$, such that $0\leq\tau^\alpha(Y_{n,\rho})-\rho<1/n$.

Let $0<\beta<\alpha$.
We may assume that  subsequences $(X_{n_k,\rho}),$ and $(Y_{n_k,\rho})$ of the sequences above converge in $C^{2,\alpha}$  and $C^{2,\beta}$, respectively,
to diffeomorphisms $X_\rho\in C^{2,1}$, $Y_\rho\in C^{2,\alpha}$ such that $\tau^1(X_\rho) \leq\mathcal{K}$,
$\tau^\alpha(Y_\rho)\leq \rho$ (cf. Lemma \ref{lbas}) and $X_\rho(T)$ is congruent to $Y_\rho(T)$. Indeed, $\tau^\alpha(Y_\rho)=\rho$. In order to prove this, we will suppose that
$\tau^\alpha(Y_\rho)<\rho$. Since  
$X_{n_k,\rho}\to X_\rho$ in $C^{2,\alpha}$ we may assume that 
there exists diffeomorphisms $Y_{n_k,\rho}^*\in  \Lambda^{X_{n_k,\rho}(T)}_{\alpha}$ converging to $Y_\rho$ in $C^{2,\alpha}$ (for instance, take $Y_{n_k,\rho}^*=X_{n_k,\rho}\circ X_\rho^{-1}\circ Y_\rho$). Hence,
by Lemma \ref{lbas} and the definition of $\Lambda^{X_{n_k,\rho}(T)}_{\alpha}$, there holds that $\rho\leq \lim_{k\to\infty}\tau^\alpha(Y_{n_k,\rho}^*)=\tau^\alpha(Y_\rho)<\rho$, a contradiction.

Let $M_\rho=Y_\rho(T)$. Then $M_\rho$ is a $C^{2,1}$ torus embedded in $\S3$. Identifying 
$X_{n,\rho}(T)$ and $Y_{n,\rho}(T)$ (that are congruents) there holds that $Y_{n,\rho}(T)$ has mean curvature equal to zero at the points
$X_{n,\rho}(p^{jk}_n)$, $1\leq j,k\leq n$. We remark that the points $X_{\rho}(p^{jk}_n)$, $n\in \mathbb N$, are dense in $M_\rho$. Thus convergence in $C^{2,\alpha}$ implies that $M_\rho$ is minimal. 

Suppose that  $M_\rho$ is congruent to the Clifford torus $T$ for some $\rho\neq 0$. Then $I(T)$ is congruent to 
$X_{\rho}(T)$. As above there exists $C^{2,\alpha}$ diffeomorphisms $Y_{n_k,\rho}^*$ converging to $I$ in $C^{2,\alpha}$ such that $Y_{n_k,\rho}^*(T)=X_{n_k,\rho}(T)$, i.e., $Y_{n_k,\rho}^*\in\Lambda^{X_{n,\rho}(T)}_{\alpha}$. A contradiction
because $\lim_{k\to\infty}\tau^\alpha(Y_{n_k,\rho}^*)=0<\rho$.
\qed

\medskip

The following corollary follows immediately from Lemma \ref {lfinal}:

\begin{corollary}\label{cfinal} Let $0<\alpha<1$. If there exists a torus $M$ minimally embedded in $\S3$  that is not congruent to the Clifford torus $T$ then there also exists a
sequence $(X_k)$  of $C^{2,\alpha}$ diffeomorphisms of $\S3$ canonically extended to $A_2$ such that each $M_k=X_k(T)$is a torus minimally embedded in $\S3$ noncongruent to  the Clifford torus
and $\lim_{k\to\infty}||X_k-I||_{C^{2,\alpha}}= 0$.
\end{corollary}

\section{Main Result} \label{sl}

Corollary \ref{cfinal} implies  
the main result in this work if it is further proved that  the Clifford torus is 
isolated in $C^{2,\alpha}$.  
In order to accomplish this we shall deal with the classical Hilbert spaces $L^2(M)$ and  $H^1(M)$ of functions along a riemannian manifold $(M;\langle\hspace{.02in},\hspace{.02in}\rangle)$. They are obtained by the completion
of $C(M)$ and $C^1(M)$ with respect to the norms derived from the respective 
inner products:
\begin{eqnarray}
\langle f,g \rangle_{L^2(M)} &=& \int_M fg\hspace{.04in} dM, \\
\langle f,g \rangle_{H^1(M)} &=&  \langle f,g\rangle_{L^2(M)} +\int_M \langle \hbox{grad}_M f,\hbox{grad}_M g\rangle\hspace{.04in} dM. 
\end{eqnarray}
We refer to \cite{Be} for the properties of these spaces that we will need below, as those regarding the eigenvalues of the Laplace-Beltrami operator $\triangle_{M}$ on $C^2$ functions along $M^2\subset\S3$, i.e.,
the nonnegative $\lambda$ for which   the equation $\triangle_{M}u+\lambda\hspace{.02in}u=0$ has a nontrivial
solution $u\in C^\infty(M)$ ($u$ is then called  an eigenfunction of $\lambda$). We recall the variational characterization of the first positive eigenvalue:

\begin{equation}\label{epa}
\lambda_1(M)=\inf_{{f\in C^1(M)}\atop{\int_M f dM=0}}\dfrac{\int_{M}|\hbox{grad}_{M} f|^2\ dM}{\int_M f^2 dM}.
\end{equation}
Finally, we remark that $M^2$ is minimally immersed in $\S3$ if and only if the coordinate functions of the 
immersion in the ambient euclidean space $\mathbb R^4$ are eigenfunctions of the Laplace-Beltrami operator in 
the induced metric with eigenvalue equal to 2 (cf. \cite{lw1}, parag. 4, prop. 1).

\begin{lemma}\label{ltiso} Let $0<\alpha<1$. Then there does not exist any sequence $(X_k)$ of diffeomorphisms of $\S3$
canonically extended to $A_2$  such that each $M_k=X_k(T)$ is a torus minimally embedded in $\S3$ noncongruent to the 
Clifford torus $T$ and $\lim_{k\to \infty}||X_k-I||_{C^{2,\alpha}}=0$.
\end{lemma}

\noindent{\it Proof.} We admit the existence of the sequences $(M_k),\ (X_k)$ described above.
Let $\lambda_1(M_k)$ be the first
positive eigenvalue of the Laplace-Beltrami operator $\triangle_{M_k}$.
From the lower bound for $\lambda_1(M_k)$ in \cite{cw} and the remark above on the coordinate functions of minimal immersions in $\S3$ there holds that  $1<\lambda_1(M_k)\leq 2$.

By taking a subsequence if necessary we may assume that $(\lambda_1(M_k))$ is a convergent sequence. It is clear
that $1\leq \lim_{k\rightarrow \infty}(\lambda_1(M_k))\leq 2$. 
Let 
$f_k:M_k\rightarrow \mathbb R^4$  be a nontrivial eigenfunction of the first
eigenvalue
$\lambda_1(M_k)$ of $M_k$ and  let $g_k=f_k\circ X_k$. Assume that
$\int_{M_k} f_k^2 \hspace{.04in}dM_k=1$. 
These assumptions  imply that $\lambda_1(M_k)=\int_{M_k} |\hbox{grad}_{M_k} f_k|^2
dM_k$ (see eq. \ref{epa}).
Let $\lambda=\lim_{k\to\infty}\lambda(M_k)$. Given $\delta>0$
the following inequality holds for $k>0$ sufficiently large
\begin{equation}\label{etp1}
0\leq\langle g_k,g_k\rangle_{H^1(T)}^2
\leq 1+\lambda +\delta,
\end{equation}
since  $\lambda_1(M_k)\rightarrow \lambda$ and the norm of the jacobian of $X_k$ tends to 1 when $k\to\infty$.
Hence the sequence $(g_k)$ is bounded in $H^{1}(T)$.
Thus we may
admit by taking a subsequence if necessary that the sequence  $(g_k)$ converges weakly to a function $g$ in $H^1(T)$.
It is well known that the inclusion $H^1(T)\subset L^2(T)$ is continuous
and
compact (see \cite{Be}). So the sequence $(g_k)$ is pre-compact in $L^2(T)$ an hence we
can find
a subsequence $(g_{k_i})$ that converges weakly to $g$ in $H^1(T)$ and
strongly to a function
$\overline g$ in $L^2(T)$. Thus $g_{k_i}$ converges weakly to $\overline
g$ in $L^2(T)$ too since strong convergence implies weak convergence.
The continuity of the inclusion  $H^1(T)\subset L^2(T)$ implies that
$g=\overline g$.
From standard regularization procedures and
the assumptions on $f_k$, $M_k$ and $g_k$, we obtain
that $g$ is a smooth function satisfying $\int_{T} g^2 \hspace{.04in}
dM_T=1$ and
$\triangle_T g +\lambda g=0$.
Thus 
$1\leq \lambda\leq 2$ is an
eigenvalue
of the Clifford torus, that is known to have the least positive eigenvalue equal to 2 (cf. \cite{MR}). Hence $\lambda=2$.

Let $x=(x_1,x_2,x_3,x_4)$ be the coordinate functions of the Clifford torus $T$. We define analogously the coordinate functions
$x^k=(x_1^k,x_2^k,x_3^k,x_4^k)$ of the minimal embedding $M_k\subset\S3$.
As was observed above, $\triangle_T\hspace{.02in} x +2x=
 \triangle_{M_k}x^k+2x^k=0$. The functions $x_i$ are a basis for the eigenspace
of the first eigenvalue of the Clifford torus (see \cite{Ch}). Thus the
function $g$
defined in the previous paragraph, an eigenfunction of the first positive eigenvalue of the Clifford torus as was shown, 
is given by $g=  \sum_{i=1}^4 a_i x_i$. Let $h_k=\sum_{i=1}^4 a_ix_i^k$.
We now observe that $\lambda_1(M_k)<2$ or $M_k$ is the Clifford torus
(Theorem 4 in \cite{MR}). Hence if $f_k$ is the nontrivial eingenfunction of the first eingenvalue $\lambda(M_k)$
described above 
then $\int_{M_k} f_kh_k\hspace{.04in} dM_k=0$. A continuity argument implies that  $\int_{T}
g^2dM_T=0$, a
contradiction. \qed
\smallskip

Theorem \ref{tl1} follows from Cor. \ref{lfinal} and Lemma \ref{ltiso}.

\begin{theorem}\label{tl1}
If $M$ is a torus minimally embedded in $\S3$ then $M$ is congruent to the
Clifford torus.
\end{theorem}

\appendix

%\newpage
\section {Intersections of equators  and tori in $\S3$}\label{afig}

We  sketch here  all  the possible intersections of a equator $S(v)$ in $\S3$ and a
embedded torus $M$ with the {\it two-piece property}. For a better understanding of these pictures
one may assume that we have performed a stereographical  projection from $\S3$
to $\R3$ with  north pole a point in $S(v)\backslash M$. Thus the equator $S(v)$ is taken
to a plane  and the image of $M$ is an ordinary torus of $\R3$. For types 3 and 4 we have assumed that the north pole
of the stereographical projection is chosen in such a way that the curves in $S(v)\cap M$ are taken to meridians of the projected torus.

\medskip
\centerline{\epsfxsize=4.0cm \epsfysize=4.0cm
 \epsfbox{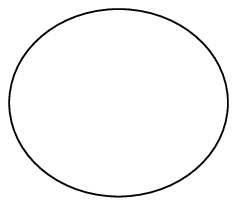}\epsfxsize=4.2cm \epsfysize=4.2cm \epsfbox{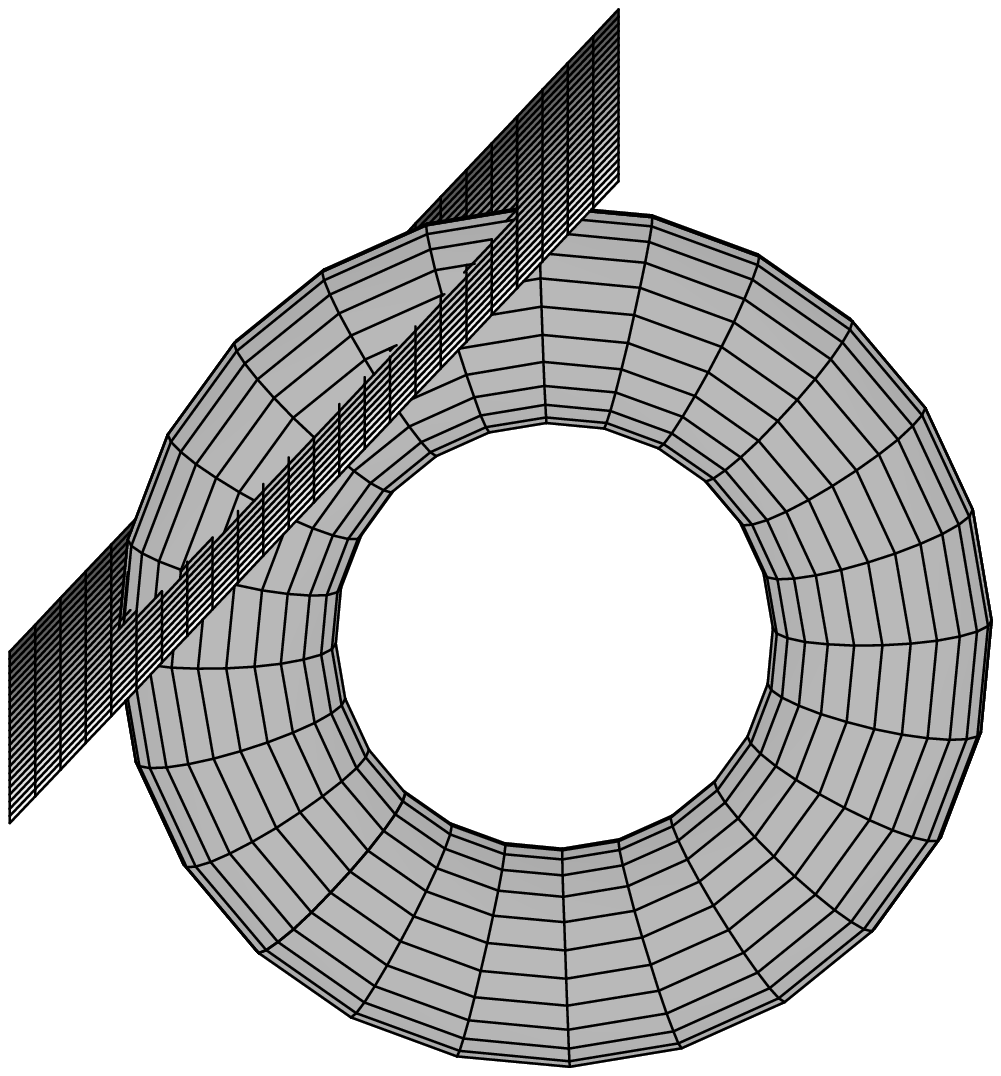}}
\smallskip
\centerline{type 1}
\medskip

\medskip
\centerline{\epsfxsize=4.0cm \epsfysize=4.0cm \epsfbox{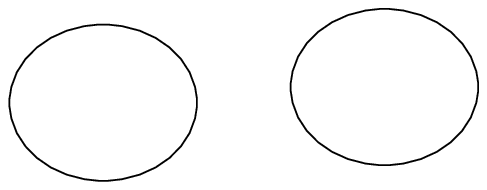}
 \epsfxsize=4.2cm \epsfysize=4.2cm \epsfbox{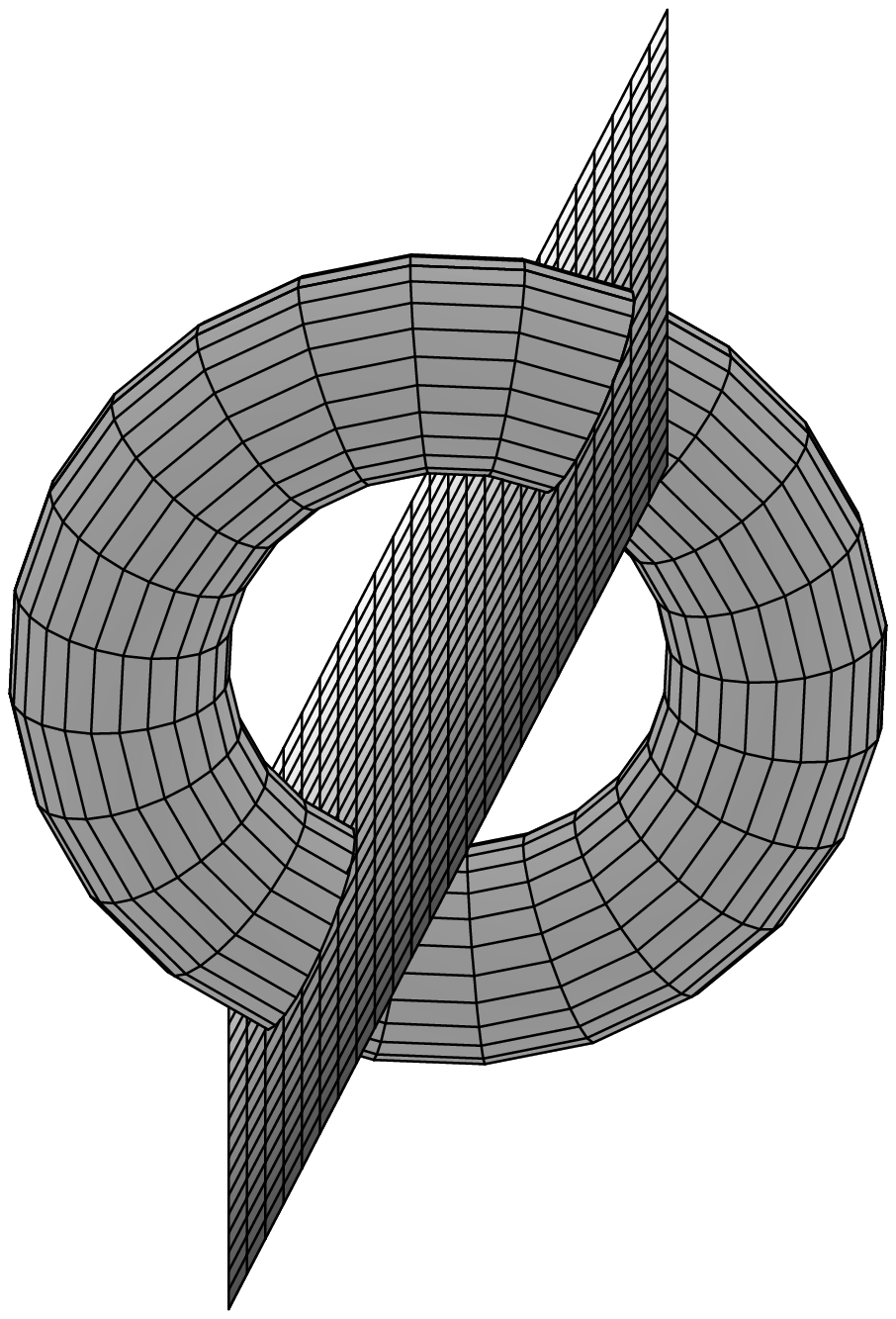}
}
\smallskip
\centerline{type 2}
\medskip

\medskip
\centerline{\epsfxsize=4.0cm \epsfysize=4.0cm
  \epsfbox{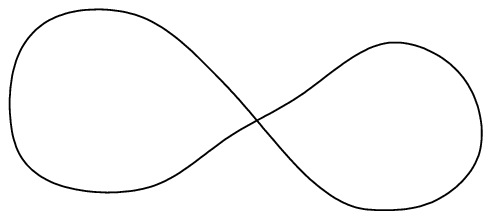}
  \epsfxsize=4.2cm \epsfysize=4.2cm \epsfbox{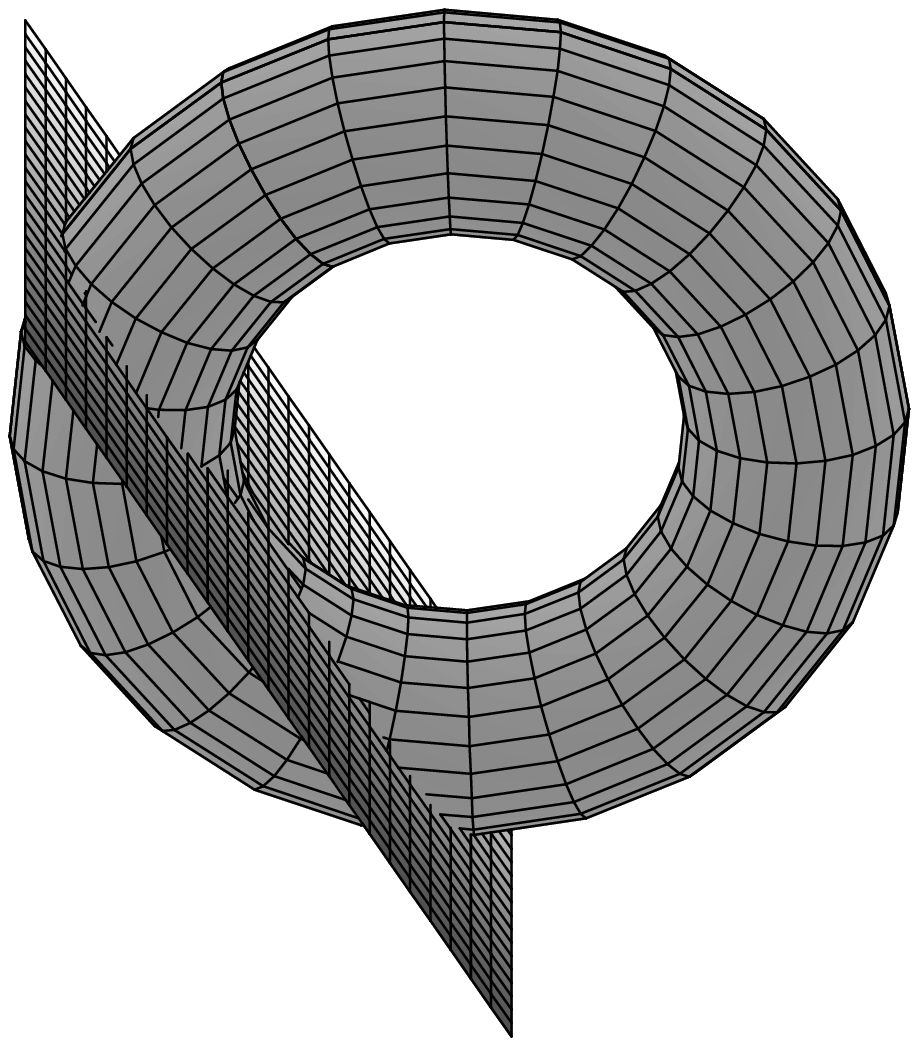}}
\smallskip
\centerline{type 3}
\medskip

\medskip
\centerline{\epsfxsize=4.0cm \epsfysize=4.0cm
  \epsfbox{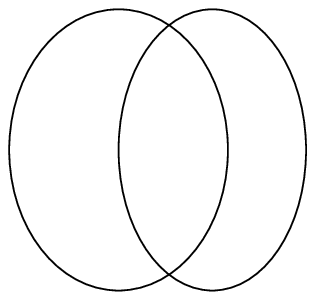}
 \epsfxsize=4.2cm \epsfysize=4.2cm \epsfbox{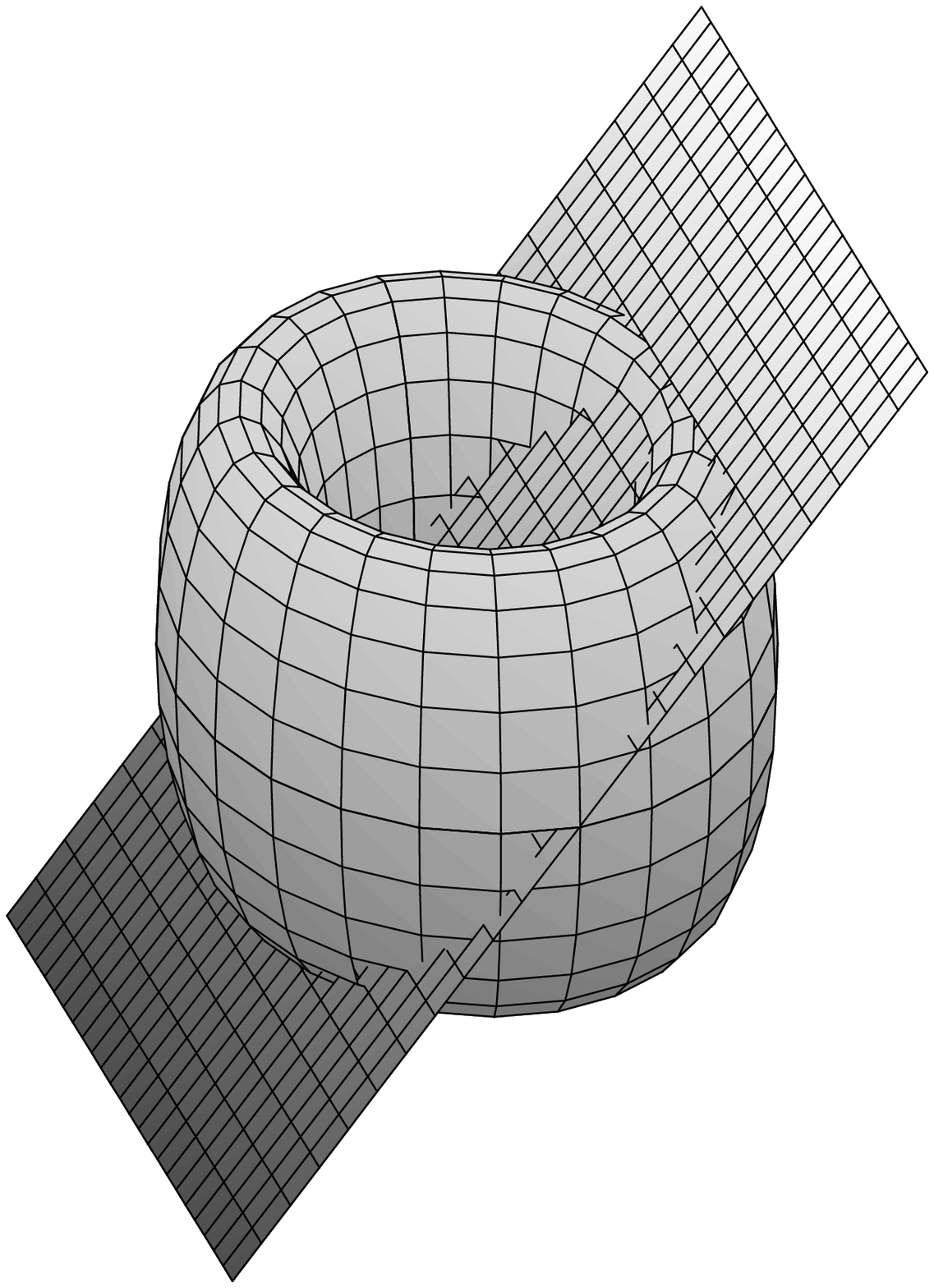}}
\smallskip
\centerline{type 4}

\medskip\bigskip

Fernando A. A. Pimentel: Departamento de Matem{\'a}tica, UFC,
Campus do Pici bloco 914,
60455-760 Fortaleza-Ce, 
Brazil \ \ \ \ \ \ \ \ \ \ {\small pimentelf@gmail.com}

\end{document}